\documentclass[aip,jmp,amsmath,amssymb,preprint]{revtex4-1}
\usepackage[dvipdfmx]{graphicx}

\begin{document}

\title{An extension of Gauss's arithmetic-geometric mean (AGM) to three variables iteration scheme}
\author{Kiyoshi Sogo}
\thanks{EMail: sogo@icfd.co.jp}
\affiliation{
Institute of Computational Fluid Dynamics, 
1-16-5, Haramachi, Meguro, Tokyo, 152-0011, Japan
}

\begin{abstract}
Gauss's arithmetic-geometric mean (AGM) which is described by two variables iteration $(a_n, b_n)\rightarrow (a_{n+1}, b_{n+1})$ by 
$a_{n+1}=(a_n+b_n)/2,\ b_{n+1}=\sqrt{a_nb_n}$. We extend it to 
three variables iteration $(a_n, b_n, c_n)\rightarrow (a_{n+1}, b_{n+1}, c_{n+1})$ which reduces to Gauss's AGM when $c_0=0$. 
Our iteration starting from $a_0>b_0>c_0>0$ with further restriction $a_0>b_0+c_0$ 
converges to $a_\infty=b_\infty=M(a_0, b_0, c_0)$ and $c_\infty=0$. 
The limit $M(a_0, b_0, c_0)$ is expressed by Appell's hyper-geometric function $F_1(1/2, \{1/2, 1/2\}, 1; \kappa, \lambda)$ 
of two variables $(\kappa, \lambda)$ which are determined by $(a_0, b_0, c_0)$. 
A relation between two hyper-geometric functions (Gauss's and Appell's) is found as a by-product.
\end{abstract}



\maketitle

\section{Introduction}
\setcounter{equation}{0}

The theory of arithmetic-geometric mean (AGM) by Gauss was published posthumously in his Werke.\cite{Gauss} 
Gauss observed that the iteration 
\begin{align}
a_{n+1}=\frac{1}{2}(a_n+b_n),\quad
b_{n+1}=\sqrt{a_n b_n},\qquad (a_0 > b_0>0)
\label{GaussAGM} 
\end{align}
gives converging series $a_0>a_1>\cdots >a_\infty=b_\infty>\cdots>b_1>b_0$.
We write the limit $a_\infty=b_\infty=\text{agM}(a_0, b_0)$ and call it arithmetic-geometric mean (AGM), which satisfies 
a functional equation
\begin{align}
\text{agM}(a, b)=\text{agM}((a+b)/2, \sqrt{ab}).
\label{GaussF1}
\end{align}
The amazing discovery \cite{Borwein} by Gauss is the equality
\begin{align}
\frac{a_0}{\text{agM}(a_0, b_0)}=\frac{1}{\pi}\int_0^1\frac{du}{\sqrt{u(1-u)(1-\kappa u)}},\quad 
\kappa=1-\left(\frac{b_0}{a_0}\right)^2
\end{align}
where the right hand side integral ($0\leq\kappa<1$) is a special case of Gauss's hyper-geometric function
\begin{align}
\frac{1}{\pi}\int_0^1\frac{du}{\sqrt{u(1-u)(1-\kappa u)}}={}_2F_1\left(\frac{1}{2}, \frac{1}{2}, 1; \kappa\right),
\label{GaussIntegral}
\end{align}
which is generally defined by
\begin{align}
{}_2F_1\left(\alpha, \beta, \gamma; \kappa\right)
&=\frac{\Gamma(\gamma)}{\Gamma(\alpha)\Gamma(\gamma-\alpha)}\int_0^1 u^{\alpha-1}(1-u)^{\gamma-\alpha-1}(1-\kappa u)^{-\beta}\ du,
\\
&=\sum_{n=0}^\infty \frac{(\alpha)_n(\beta)_n}{(\gamma)_n}\cdot\frac{\kappa^n}{n!},
\qquad (|\kappa|<1)
\end{align}
where $\gamma>\alpha>0$ is assumed, and Pochhammer's symbol $(\alpha)_0=1,\ (\alpha)_1=\alpha,\ (\alpha)_n=\alpha(\alpha+1)\cdots(\alpha+n-1)=\Gamma(\alpha+n)/\Gamma(\alpha)$ is used. To derive this series expansion, 
the following Newton's expansion formula is used
\begin{align}
(1-\kappa u)^{-\beta}=\sum_{m=0}^\infty \left(\frac{(\beta)_n}{n!}\ \kappa^n\right)\cdot u^n.
\label{Newton}
\end{align}

To extend the above results, we start by considering the integral with an additional new parameter $\lambda$
\begin{align}
\frac{1}{\pi}\int_0^1 \frac{du}{\sqrt{u(1-u)(1-\kappa u)(1-\lambda u)}},\qquad (0<\lambda<\kappa<1)
\label{newGaussIntegral}
\end{align}
which reduces to \eqref{GaussIntegral} when $\lambda=0$. In the next section we construct a model with the iteration rule 
$(a_n, b_n, c_n)\rightarrow (a_{n+1}, b_{n+1}, c_{n+1})$ and show the converging limit $M(a_0, b_0, c_0)$ 
is given by
\begin{align}
\frac{a_0}{M(a_0, b_0, c_0)}=\frac{1}{\pi}\int_0^1 \frac{du}{\sqrt{u(1-u)(1-\kappa u)(1-\lambda u)}},
\label{Result}
\end{align}
where the relations of $(\kappa, \lambda)$ with $(a_0, b_0, c_0)$ are given by
\begin{align}
\left(\frac{b_0}{a_0}\right)^2=(1-\kappa)(1-\lambda),\qquad \left(\frac{c_0}{a_0}\right)^2=\kappa\lambda,
\end{align}
which will be derived in the next section.

\section{Extension of Gauss's AGM}
\setcounter{equation}{0}

\subsection{Iteration rule}

Let us begin to observe that the integral \eqref{GaussIntegral} of Gauss is rewritten, by setting $u=\sin^2\theta$, 
\begin{align}
\frac{2a_0}{\pi}\int_0^{\pi/2}\frac{d\theta}{\sqrt{a_0^2\cos^2\theta+b_0^2\sin^2\theta}}
=\frac{1}{\pi}\int_0^1\frac{du}{\sqrt{u(1-u)(1-\kappa u)}},\quad \kappa=1-\left(\frac{b_0}{a_0}\right)^2,
\end{align}
which can be extended to
\begin{align}
&\frac{2a_0}{\pi}\int_0^{\pi/2}\frac{d\theta}{\sqrt{a_0^2\cos^2\theta+b_0^2\sin^2\theta-c_0^2\cos^2\theta\sin^2\theta}}
\nonumber \\
&\qquad
=\frac{1}{\pi}\int_0^1\frac{du}{\sqrt{u(1-u)(1-\kappa u)(1-\lambda u)}},
\end{align}
by introducing a new variable $c_0$.
Here we have, by substituting $\sin^2\theta=u$,
\begin{align}
&a_0^2\cos^2\theta+b_0^2\sin^2\theta-c_0^2\cos^2\theta\sin^2\theta=a_0^2-(a_0^2-b_0^2+c_0^2)u+c_0^2u^2
\nonumber \\
&\qquad =a_0^2(1-\kappa u)(1-\lambda u)
\quad\Longrightarrow\quad 
(1-\kappa)(1-\lambda)=\left(\frac{b_0}{a_0}\right)^2,\quad \kappa\lambda=\left(\frac{c_0}{a_0}\right)^2,
\end{align}
from which we can express $\kappa, \lambda$ in terms of $a_0, b_0, c_0$ as follows
\begin{align}
\begin{split}
\kappa=\frac{1}{2}\left\{1-\frac{b_0^2-c_0^2}{a_0^2}+
\sqrt{\left(1-\frac{(b_0+c_0)^2}{a_0^2}\right)\left(1-\frac{(b_0-c_0)^2}{a_0^2}\right)}\ \right\},
\\
\lambda=\frac{1}{2}\left\{1-\frac{b_0^2-c_0^2}{a_0^2}-
\sqrt{\left(1-\frac{(b_0+c_0)^2}{a_0^2}\right)\left(1-\frac{(b_0-c_0)^2}{a_0^2}\right)}\ \right\},
\end{split}
\label{KL1}
\end{align}
and 
\begin{align}
\begin{split}
1-\kappa=\frac{1}{2}\left\{1+\frac{b_0^2-c_0^2}{a_0^2}-
\sqrt{\left(1-\frac{(b_0+c_0)^2}{a_0^2}\right)\left(1-\frac{(b_0-c_0)^2}{a_0^2}\right)}\ \right\},
\\
1-\lambda=\frac{1}{2}\left\{1+\frac{b_0^2-c_0^2}{a_0^2}+
\sqrt{\left(1-\frac{(b_0+c_0)^2}{a_0^2}\right)\left(1-\frac{(b_0-c_0)^2}{a_0^2}\right)}\ \right\},
\end{split}
\label{KL2}
\end{align}
where conditions $a_0>b_0>c_0>0$ and $a_0>b_0+c_0$ are assumed to guarantee $\kappa$ and $\lambda$ take real values
($0<\lambda<\kappa<1$).
If we introduce $\xi, \eta$ by
\begin{align}
\xi=\frac{b_0+c_0}{a_0},\qquad
\eta=\frac{b_0-c_0}{a_0},
\end{align}
eqs. \eqref{KL1} and \eqref{KL2} are written in much simpler forms by
\begin{align}
\begin{split}
\kappa&=\frac{1}{2}\left(1-\xi\eta+\sqrt{(1-\xi^2)(1-\eta^2)}\ \right),\\ 
\lambda&=\frac{1}{2}\left(1-\xi\eta-\sqrt{(1-\xi^2)(1-\eta^2)}\ \right),
\end{split}
\\
\begin{split}
1-\kappa&=\frac{1}{2}\left(1+\xi\eta-\sqrt{(1-\xi^2)(1-\eta^2)}\ \right),\\ 
1-\lambda&=\frac{1}{2}\left(1+\xi\eta+\sqrt{(1-\xi^2)(1-\eta^2)}\ \right).
\end{split}
\end{align}
Therefore $\kappa,\ \lambda$, which are $\kappa_0, \lambda_0$ actually, 
are determined by $(a_0, b_0, c_0)$ through $\xi, \eta$. 

To construct the iteration rule, we notice, when $0\leq\lambda<\kappa<1$,
\begin{align}
\frac{1}{\pi}\int_0^1\frac{du}{\sqrt{u(1-u)(1-\kappa u)(1-\lambda u)}}=
\frac{2}{\pi\sqrt{1-\lambda}}\cdot K(k_0),\qquad k_0^2=1-\frac{1-\kappa}{1-\lambda}
\label{GaussAppell}
\end{align}
where $K(k)$ is the complete elliptic integral of the first kind
\begin{align}
K(k)=\int_0^{\pi/2}\frac{d\phi}{\sqrt{1-k^2\sin^2\phi}}=\frac{1}{2}\int_0^1\frac{du}{\sqrt{u(1-u)(1-k^2 u)}}.
\end{align}

We use here the Landen transformation formula \cite{Borwein}
\begin{align}
K(k_0)=(1+k_1)\cdot K(k_1),
\qquad k_1=\frac{1-\sqrt{1-k_0^2}}{1+\sqrt{1-k_0^2}},
\end{align}
which can be rewritten such as
\begin{align}
K(k_1)=\frac{1+\sqrt{1-k_0^2}}{2}\cdot K(k_0),
\end{align}
because we have
\begin{align}
1+k_1=\frac{2}{1+\sqrt{1-k_0^2}}.
\label{k0k1}
\end{align}
Therefore our iteration rule, $(\kappa_0, \lambda_0)\rightarrow (\kappa_1, \lambda_1)$, is
\begin{align}
&\frac{1}{a_0\pi}\int_0^1\frac{du}{\sqrt{u(1-u)(1-\kappa_0 u)(1-\lambda_0 u)}}=
\frac{1}{a_1\pi}\int_0^1\frac{du}{\sqrt{u(1-u)(1-\kappa_1 u)(1-\lambda_1 u)}}
\nonumber \\
&\Longrightarrow\quad
\frac{1}{a_0\sqrt{1-\lambda_0}}\cdot K\left(k_0=\sqrt{1-\frac{1-\kappa_0}{1-\lambda_0}}\ \right)=
\frac{1}{a_1\sqrt{1-\lambda_1}}\cdot K\left(k_1=\sqrt{1-\frac{1-\kappa_1}{1-\lambda_1}}\ \right)
\nonumber \\
&\Longrightarrow\quad
\frac{1}{a_0\sqrt{1-\lambda_0}}\cdot\frac{2}{1+\sqrt{1-k_0^2}}=\frac{1}{a_1\sqrt{1-\lambda_1}},
\end{align}
which is rewritten into
\begin{align}
\frac{a_1}{a_0}&=\frac{1+\sqrt{1-k_0^2}}{2}\cdot\sqrt{\frac{1-\lambda_0}{1-\lambda_1}},\qquad
k_0^2=1-\frac{1-\kappa_0}{1-\lambda_0},
\nonumber \\
&=\frac{\sqrt{1-\lambda_0}+\sqrt{1-\kappa_0}}{2 \sqrt{1-\lambda_1}}.
\end{align}

Now it is convenient to introduce parameters $\alpha$ and $\varepsilon$ ($0<\varepsilon<\alpha<\pi/2$) by
\begin{align}
\kappa=\sin^2\alpha,\qquad \lambda=\sin^2\varepsilon,
\end{align}
which imply
\begin{align}
\begin{split}
&\frac{b}{a}=\cos\alpha\cos\varepsilon,\quad 
\frac{c}{a}=\sin\alpha\sin\varepsilon,\\
&\frac{b+c}{a}=\cos(\alpha-\varepsilon),\quad
\frac{b-c}{a}=\cos(\alpha+\varepsilon).
\end{split}
\end{align}
By use of these parameters, we have
\begin{align}
k_0^2=1-\frac{1-\kappa_0}{1-\lambda_0}=1-\frac{\cos^2\alpha_0}{\cos^2\varepsilon_0}=
\frac{\sin(\alpha_0+\varepsilon_0)\sin(\alpha_0-\varepsilon_0)}{\cos^2\varepsilon_0},
\end{align}
and similarly
\begin{align}
k_1^2=1-\frac{1-\kappa_1}{1-\lambda_1}=
\frac{\sin(\alpha_1+\varepsilon_1)\sin(\alpha_1-\varepsilon_1)}{\cos^2\varepsilon_1}.
\end{align}
Substituting these into \eqref{k0k1}, we have
\begin{align}
\begin{split}
&\tan^2\left(\frac{\alpha_0+\varepsilon_0}{2}\right)=\frac{\sin(\alpha_1+\varepsilon_1)}{\cos\varepsilon_1}=
\sin\alpha_1+\cos\alpha_1\tan\varepsilon_1,
\\
&\tan^2\left(\frac{\alpha_0-\varepsilon_0}{2}\right)=\frac{\sin(\alpha_1-\varepsilon_1)}{\cos\varepsilon_1}=
\sin\alpha_1-\cos\alpha_1\tan\varepsilon_1,
\end{split}
\end{align}
from which we obtain
\begin{align}
\begin{split}
&\sin\alpha_1=\frac{1}{2}\left\{\tan^2\left(\frac{\alpha_0+\varepsilon_0}{2}\right)+
\tan^2\left(\frac{\alpha_0-\varepsilon_0}{2}\right)\right\},
\\
&\tan^2\varepsilon_1=\tan^2\alpha_1-\sec^2\alpha_1\cdot\tan^2\left(\frac{\alpha_0+\varepsilon_0}{2}\right)\cdot\tan^2\left(\frac{\alpha_0-\varepsilon_0}{2}\right),
\end{split}
\end{align}
which gives $\kappa_1=\sin^2\alpha_1$ and $\lambda_1=\sin^2\varepsilon_1$ 
in terms of $(\alpha_0, \varepsilon_0)$. 
These are the iteration scheme what we are looking for. 

Here we can rewrite these quantities in terms of $\xi, \eta$. Since 
\begin{align}
\cos(\alpha_0-\varepsilon_0)=\frac{b_0+c_0}{a_0}=\xi,\quad
\cos(\alpha_0+\varepsilon_0)=\frac{b_0-c_0}{a_0}=\eta,
\end{align}
we have at first
\begin{align}
&\tan^2\left(\frac{\alpha_0-\varepsilon_0}{2}\right)=\frac{1-\xi}{1+\xi},\quad
\tan^2\left(\frac{\alpha_0+\varepsilon_0}{2}\right)=\frac{1-\eta}{1+\eta}
\\
&\quad\Longrightarrow\quad
\sin\alpha_1=\frac{1}{2}\left(\frac{1-\xi}{1+\xi}+\frac{1-\eta}{1+\eta}\right)=\frac{1-\xi\eta}{(1+\xi)(1+\eta)},
\end{align}
therefore we obtain
\begin{align}
\kappa_1=\sin^2\alpha_1=\left(\frac{1-\xi\eta}{(1+\xi)(1+\eta)}\right)^2.
\end{align}
We have next, after some calculations 
\begin{align}
\tan^2\varepsilon_1&=\frac{\sin^2\alpha_1}{1-\sin^2\alpha_1}-\frac{1}{1-\sin^2\alpha_1}\cdot\frac{(1-\xi)(1-\eta)}{(1+\xi)(1+\eta)}
\nonumber \\
&=\frac{(\xi-\eta)^2}{(2+\xi+\eta)(\xi+\eta+2\xi\eta)},
\end{align}
therefore we obtain
\begin{align}
\lambda_1=\sin^2\varepsilon_1=\frac{\tan^2\varepsilon_1}{1+\tan^2\varepsilon_1}=\frac{(\xi-\eta)^2}{2(1+\xi)(1+\eta)(\xi+\eta)}.
\end{align}

In conclusion, the procedure of iteration is described as follows.
\begin{align}
&\quad (1)\ \text{Compute}\ \kappa_0,\ \lambda_0\ \text{by using}\ \xi=(b_0+c_0)/a_0,\ \eta=(b_0-c_0)/a_0,\ 
\nonumber \\
\begin{split}
&\qquad\qquad \kappa_0=\frac{1}{2}\left(1-\xi\eta+\sqrt{(1-\xi^2)(1-\eta^2)}\right),
\\
&\qquad\qquad \lambda_0=\frac{1}{2}\left(1-\xi\eta-\sqrt{(1-\xi^2)(1-\eta^2)}\right).
\end{split}
\\
&\quad (2)\ \text{Compute}\ \kappa_1,\lambda_1\ \text{by using the same}\ \xi,\ \eta,
\nonumber \\
\begin{split}
&\qquad\qquad \kappa_1=\left(\frac{1-\xi\eta}{(1+\xi)(1+\eta)}\right)^2,
\\
&\qquad\qquad \lambda_1=\frac{(\xi-\eta)^2}{2(1+\xi)(1+\eta)(\xi+\eta)}.
\end{split}
\\
&\quad (3)\ \text{Compute}\ (a_1, b_1, c_1)\ \text{by using}\ \kappa_0, \lambda_0\ \text{and}\ \kappa_1, \lambda_1
\nonumber \\
\begin{split}
&\qquad\qquad a_1=a_0\ (\sqrt{1-\lambda_0}+\sqrt{1-\kappa_0})/2\sqrt{1-\lambda_1},
\\
&\qquad\qquad b_1=a_1\ \sqrt{(1-\kappa_1)(1-\lambda_1)},
\\
&\qquad\qquad c_1=a_1\ \sqrt{\kappa_1\lambda_1}.
\end{split}
\end{align}
We perform these sequences repeatedly, regarding the sub-indices $(0, 1)$ as $(n, n+1)$. 
The result of iterations is $a_\infty=b_\infty=M(a_0, b_0, c_0)$ and $c_\infty=0$ as was announced before.

Let us give an example of numerical computation by Python, whose code is given in Appendix. 
For the case of $a_0=1.0,\ b_0=0.5,\ c_0=0.2$, after ten iterations we get 
$a_{10}=b_{10}=0.7250921711406717$ and $c_{10}=0.0$, as is expected. 
On the other hand, the numerical integration of \eqref{Result} gives $M(a_0, b_0, c_0)=0.725092$, which can be 
regarded as the same value.

\subsection{Appell's hyper-geometric function}

Our integration \eqref{newGaussIntegral} is expressed by
\begin{align}
\frac{1}{\pi}\int_0^1\frac{du}{\sqrt{u(1-u)(1-\kappa u)(1-\lambda u)}}=
F_1\left(\frac{1}{2},\{\frac{1}{2}, \frac{1}{2}\}, 1; \kappa, \lambda\right),
\end{align}
which is a special case of Appell's hyper-geometric function of two variables defined by \cite{WhittakerWatson}
\begin{align}
F_1\left(\alpha,\{\beta, \beta'\}, \gamma; x, y\right)&=
\frac{\Gamma(\gamma)}{\Gamma(\alpha)\Gamma(\gamma-\alpha)}\int_0^1 u^{\alpha-1}(1-u)^{\gamma-\alpha-1}
(1-x u)^{-\beta}(1-y u)^{-\beta'}\ du
\nonumber \\
&=\sum_{m=0}^\infty \sum_{n=0}^\infty \frac{(\alpha)_{m+n}(\beta)_m(\beta')_n}{(\gamma)_{m+n}} \cdot\frac{x^m y^n}{m!n!},
\end{align}
where is used Newton's expansion formula \eqref{Newton} twice.
Appell's function $F_1$ satisfies the differential equations \cite{WhittakerWatson}
\begin{align}
\begin{split}
&x(1-x)\frac{\partial^2F_1}{\partial x^2}+y(1-x)\frac{\partial^2F_1}{\partial x\partial y}+
[\gamma-(\alpha+\beta+1)x]\frac{\partial F_1}{\partial x} -\beta y\frac{\partial F_1}{\partial y}-\alpha\beta F_1=0,
\\
&y(1-y)\frac{\partial^2F_1}{\partial y^2}+x(1-y)\frac{\partial^2F_1}{\partial y\partial x}+
[\gamma-(\alpha+\beta'+1)y]\frac{\partial F_1}{\partial y} -\beta' x\frac{\partial F_1}{\partial x}-\alpha\beta' F_1=0.
\end{split}
\end{align}

From the equality \eqref{GaussAppell} we have an interesting relation
\begin{align}
F_1\left(\frac{1}{2},\{ \frac{1}{2}, \frac{1}{2}\}, 1; x, y\right)=\quad 
\begin{split}
&\frac{1}{\sqrt{1-y}}\cdot {}_2F_1\left(\frac{1}{2}, \frac{1}{2}, 1; \frac{x-y}{1-y}\right),\quad (0<y \leq x<1)
\\
&\frac{1}{\sqrt{1-x}}\cdot {}_2F_1\left(\frac{1}{2}, \frac{1}{2}, 1; \frac{y-x}{1-x}\right),\quad (0<x \leq y<1)
\end{split}
\end{align}
between Appell's and Gauss's hyper-geometric functions. 
Especially when we set $y=x$, the right hand side becomes $1/\sqrt{1-x}$, which is verified by calculating 
the left hand side integral. 

In general such kind of relationship can be derived at least under the condition $\gamma=\beta+\beta'$, 
which is satisfied in our case ($1=1/2+1/2$).   
Such formula is given by
\begin{align}
F_1\left(\alpha,\{\beta, \beta'\}, \gamma; x, y\right)=\quad 
\begin{split}
(1-y)^{-\alpha}\cdot {}_2F_1\left(\alpha, \beta, \gamma; \frac{x-y}{1-y}\right),\quad
(0<y \leq x<1)
\\
(1-x)^{-\alpha}\cdot {}_2F_1\left(\alpha, \beta', \gamma; \frac{y-x}{1-x}\right),\quad
(0<x \leq y<1)
\end{split}
\label{G=A}
\end{align}
which is derived for $0<y \leq x<1$ case by the change of integration variable $u\rightarrow v=(1-y)u/(1-yu)$, under 
the condition $\gamma=\beta+\beta'$. 

It should be noted here that general transformation formulas are known \cite{Koornwinder}
for much wider class of multi-variable hyper-geometric functions, 
where the above relation \eqref{G=A} is a special case of such formulas. 

\section{Summary}
\setcounter{equation}{0}

The arithmetic-geometric mean by Gauss is extended to three variables iteration 
$(a_n, b_n, c_n)\rightarrow (a_{n+1}, b_{n+1}, c_{n+1})$, with $a_\infty=b_\infty=M(a_0, b_0, c_0),\ c_\infty=0$. 
Our main result is expressed by the equalities such as
\begin{align}
\frac{a_0}{M(a_0, b_0, c_0)}&=
\frac{2a_0}{\pi}\int_0^{\pi/2}\frac{d\theta}{\sqrt{a_0^2\cos^2\theta+b_0^2\sin^2\theta-c_0^2\cos^2\theta\sin^2\theta}},
\\
&=\frac{1}{\pi}\int_0^1\frac{du}{\sqrt{u(1-u)(1-\kappa u)(1-\lambda u)}},
\\
&=F_1\left(\frac{1}{2},\{ \frac{1}{2}, \frac{1}{2}\}, 1; \kappa, \lambda\right)
\\
&=\frac{1}{\sqrt{1-\lambda}}\cdot 
{}_2F_1\left(\frac{1}{2}, \frac{1}{2}, 1; \frac{\kappa-\lambda}{1-\lambda}\right),\quad (0 < \lambda<\kappa<1)
\end{align}
where the relations among  $(\kappa, \lambda)$ and $(a_0, b_0, c_0)$ are given by
\begin{align}
(1-\kappa)(1-\lambda)=\left(\frac{b_0}{a_0}\right)^2,\quad
\kappa\lambda=\left(\frac{c_0}{a_0}\right)^2.
\end{align}
It should be noted that the conditions $a_0>b_0>c_0>0$ and $a_0>b_0+c_0$ are assumed.

\newpage

\noindent
{\bf APPENDIX}

\begin{verbatim}
# Python code of extended Gauss-agM
import math
# initial variables
a0, b0, c0=1.0, 0.5, 0.2
# iterations
for n in range(0, 10):
    print (a0, b0, c0)
    x, y=(b0+c0)/a0, (b0-c0)/a0
    dum1=1-x*y
    dum2=math.sqrt((1-x*x)*(1-y*y))
    K0=(dum1+dum2)/2
    L0=(dum1-dum2)/2
#
    dum1=1-x*y
    dum2=(1+x)*(1+y)
    K1=math.pow(dum1/dum2, 2)
    dum1=(x-y)*(x-y)
    dum2=2*(1+x)*(1+y)*(x+y)
    L1=dum1/dum2
#
    dum1=math.sqrt(1-L0)+math.sqrt(1-K0)
    dum2=2*math.sqrt(1-L1)
    a1=a0*dum1/dum2
    b1=a1*math.sqrt((1-K1)*(1-L1))
    c1=a1*math.sqrt(K1*L1)
#
    a0, b0, c0=a1, b1, c1
# end of code
\end{verbatim}

\end{document}